\documentclass[11pt,reqno]{amsart}
\usepackage[utf8x]{inputenc}

\usepackage{hyperref}
\hypersetup{
    colorlinks=true,       
    linkcolor=blue,          
    citecolor=blue,        
    filecolor=blue,      
    urlcolor=cyan  
}
\usepackage{latexsym,amssymb,amsmath,amsthm,bbm}
\usepackage{epsfig,mathrsfs}
\usepackage{mathtools}
\usepackage{mathabx}

\newcommand{\R}{{\ensuremath{\mathbb{R}}}}
\newcommand{\N}{{\ensuremath{\mathbb{N}}}}

\renewcommand{\P}{\ensuremath{\mathbb{P}}}
\renewcommand{\dj}{d\kern-0.4em\char"16\kern-0.1em}
\newcommand{\E}{\ensuremath{\mathbb{E}}}

\newcommand{\wh}[1]{{\widehat{#1}}}

\newcommand{\sig}{\ensuremath{\mathcal{F}}}

\newcommand{\cm}{\ensuremath{\mathcal{CM}}}
\newcommand{\bbf}{\ensuremath{\mathcal{BF}}}
\newcommand{\cbf}{\ensuremath{\mathcal{CBF}}}

\newtheorem{Thm}{Theorem}[section]
\newtheorem{Cor}[Thm]{Corollary}
\newtheorem{Lem}[Thm]{Lemma}
\newtheorem{Prop}[Thm]{Proposition}

\theoremstyle{definition}
\newtheorem{Rem}[Thm]{Remark}

\theoremstyle{definition}
\newtheorem{Ex}[Thm]{Example}

\theoremstyle{definition}

\theoremstyle{definition}

\begin{document}
\numberwithin{equation}{section}
\bibliographystyle{amsalpha}

\title[Exponential decay of measures and Tauberian theorems]{Exponential decay of measures and Tauberian theorems}
\begin{abstract}
We study behavior of a measure on $[0,\infty)$ by considering its Laplace transform. If it is possible to extend the Laplace transform to a complex half-plane containing the imaginary axis, then the exponential decay of the tail of the measure occurs and under certain assumptions we show that the rate of the decay is given by the so called abscissa of convergence and extend the result of Nakagawa from \cite{Na}.  Under stronger assumptions we give behavior of density of the measure by considering its Laplace transform. 
In situations when there is no exponential decay we study occurrence of heavy tails and give an application in the theory of non-local equations. 
\end{abstract}

\author{Ante Mimica}
\address{University of Zagreb, Department of Mathematics, Bijeni\v cka cesta 30, 10000 Zagreb, Croatia}
\curraddr{}
\thanks{This work has been supported in part by Croatian Science Foundation under the project 3526}
\email{amimica@math.hr}



\maketitle

\allowdisplaybreaks[3]

\section{Introduction}

Theorems that give behavior of a positive measure $\nu$ defined on $[0,\infty)$ by using behavior of its Laplace transform defined by  $\mathcal{L}\nu(\lambda):=\int_{[0,\infty)}e^{-\lambda t}\nu(dt)\,,\ \lambda>0$ are called Tauberian theorems. A classical example is Karamata's Tauberian theorem saying that regular variation of the Laplace transform (at the origin) implies regular variation of density (or tail) of the measure (at infinity) (see \cite[Section XIII.5]{Fel},\cite{BGT}). 

Tauberian theorems concerning exponential growth also exist.  One of the first such theorems is Ikehara's Tauberian theorem (\cite[V.17,Theorem 17]{Widd}). 
\begin{Thm}[Ikehara]
 Let $\varphi\colon [0,\infty)\rightarrow [0,\infty)$ be a non-decreasing function such that 
 \[
    f(z)=\int\limits_0^\infty e^{-z t}\varphi(t)\,dt\quad  \text{ is finite for }\quad z=x+iy, \ x>1
 \]
 and the limit $\lim\limits_{x\to 1+}\left[f(x+iy)-\frac{1}{x+iy-1}\right]$ exists uniformly on every interval $-a\leq y\leq a$, $a>0$\,.
 Then
 \[
  \lim_{t\to\infty}e^{-t}\varphi(t)=1\,.
 \]
\end{Thm}
Ikehara's theorem can be rephrased as follows.  If the Laplace transform $f$ of the function $\varphi$ is analytic in the complex half-plane $\mathrm{Re}\,z>1$ having a pole at $z_0=1$ with residue 1, then $\varphi(t)\sim e^t$ as $t\to\infty$\,. One can say that Tauberian nature in Ikehara theorem is given through the singularity of the Laplace transform. There are certain generalizations of this theorem (e.g. Graham-Vaaler's theorem) which belong to the class of complex Tauberian theorems (see \cite[Chapter III]{Kor}).

While Ikehara's theorem and its generalizations deal with exponential growth, the aim of this article is to investigate Tauberian nature concerning exponential decay. 

Before stating our main results let us introduce concepts of the abscissa of convergence and a completely monotone function.
If $\nu$ is a measure on $[0,\infty)$, then it is known that there exists  $\sigma_0=\sigma_0^f\in [-\infty,\infty)$ such that the integral 
\[
f(z)=\int\limits_{[0,\infty)}e^{-z t}\,\nu(dt)
\]
converges for $\mathrm{Re}\,z>\sigma_0$, diverges for $\mathrm{Re}\,z<\sigma_0$ and has a singularity at $\sigma_0$ (see \cite[p. 37 and p. 58]{Widd}). The number $\sigma_0$ is known as the abscissa of convergence. Furthermore, $f$ is analytic in the half-plane $\mathrm{Re}\,z>\sigma_0$ and 
\begin{equation}\label{eq:extension}
	(-1)^n f^{(n)}(\sigma_0+\lambda)=\int_{[0,\infty)}e^{-(\sigma_0+\lambda) t} t^n \nu(dt)\quad \text{ for all }\quad \lambda>0\,,
\end{equation}
where $f^{(n)}$ denotes the $n$-th derivative of $f$\,.
In Ikehara's theorem,  $\sigma_0=1$ and this was the rate of the exponential growth of the corresponding function. The abscissa of convergence will continue to play a similar role in our main result. 

Laplace transforms of measures on $[0,\infty)$ belong to the class of completely monotone functions. This is the following class of functions  \[\cm:=\{f:(0,\infty)\rightarrow (0,\infty): f \text{ is a } C^\infty\text{-function and } (-1)^nf^{(n)}\geq 0\text{ for all } n\in \N\}\,.\]
Converesely, if $f\in \cm$, there exists a unique measure $\nu$ on $[0,\infty)$, called the representing measure of $f$,  such that $f=\mathcal{L}\nu$. This result is known as Bernstein's theorem (see  \cite[Theorem XIII.4.1]{Fel},\cite[Theorem 1.4]{SSV}). Since every $f\in \cm$ can be understood as the Laplace transform of a measure, we can also talk about the abscissa of convergence of $f$, analytic extension of $f$ to the set $(\sigma_0,\infty)$ and (\ref{eq:extension}) continues to hold. 

In our first result we are going to see that the abscissa of convergence continues to play the same type of role as in the Ikehara's theorem; it determines the rate of decay of the tail of the representing measure. 
\begin{Thm}\label{thm:tail}
Let $f\in \cm$ with the abscissa of convergence $\sigma_0\in (-\infty,0]$ and the representing measure $\nu$. Assume that there exists $n\in \N_0$ satisfying
\begin{equation}\label{eq:cond-tail}
	\limsup_{\lambda\to0+}\lambda\log\left[(-1)^n f^{(n)}(\sigma_0+\lambda)\right]=0\quad \text{ and }\quad \limsup_{\lambda\to0+}\frac{(-1)^nf^{(n)}(\sigma_0+2\lambda)}{(-1)^nf^{(n)}(\sigma_0+\lambda)}<1\,.
\end{equation}
Then
\[
	\lim_{t\to\infty}\frac{1}{t}\log\nu(t,\infty)=\sigma_0\,.
\]
\end{Thm}

In other words, Theorem \ref{thm:tail} says that if the Laplace transform of a measure can be analytically extended beyond the imaginary axis in the complex plane, under some mild conditions, exponential decay of the tail of the measure occurs. The second condition in (\ref{eq:cond-tail}) is equivalent to  the following condition (see Lemma \ref{lem:equiv}):

there exist $c>0$, $\gamma>0$ and $\Lambda_0>0$ such that 
	\begin{equation}\label{eq:upper_scaling}
		\frac{(-1)^nf^{(n)}(\sigma_0+\lambda_2)}{(-1)^nf^{(n)}(\sigma_0+\lambda_1)}\leq c\left(\frac{\lambda_2}{\lambda_1}\right)^{-\gamma}\text{ for all } 0<\lambda_1\leq \lambda_2\leq \Lambda_0\,.
	\end{equation}
Condition (\ref{eq:upper_scaling}) can be understood as a variant of upper scaling of the function $\lambda\mapsto (-1)^n f^{(n)}(\sigma_0+\lambda)$. The first condition in (\ref{eq:cond-tail}) is to ensure that this function does not explode exponentially as $\lambda$ goes to $0$\,.

The role of derivatives in  (\ref{eq:cond-tail}) is important because in some cases it begins to hold only when we start to take derivatives. It turns out that this condition is quite general and it holds in many situations, e.g.  in the case of regular variation. Recall that a function $h\colon (0,\infty)\rightarrow (0,\infty)$ varies regularly at $0$ 
with index $\rho\in\R$ if 
\[
	\lim_{\lambda\to 0+}\frac{h(\lambda x)}{h(x)}=x^{\rho}\quad \text{ for all }\quad x>0\,.
\]
If $h\colon (0,\infty)\rightarrow (0,\infty)$ varies regularly with index $0$, we say that it varies slowly. It is known that for any regularly varying function $h\colon (0,\infty)\rightarrow (0,\infty)$  with index $\rho\in \R$ there exists a slowly varying function $\ell\colon (0,\infty)\rightarrow (0,\infty)$ such that
\begin{equation}\label{eq:rv-repr}
	h(\lambda)=\lambda^\rho\ell(\lambda)\quad\text{ for all }\quad  \lambda>0
\end{equation}
(cf. \cite[Theorem 1.4.1]{BGT})\,.

\begin{Cor}\label{cor:tail_reg-var}
Let $f\in \cm$ with the representing measure $\nu$ and the abscissa of convergence $\sigma_0\in (-\infty,0]$ and assume that there exists $n\in \N_0$ such that the function $h\colon (0,\infty)\rightarrow(0,\infty)$ defined by $h(\lambda):=(-1)^nf^{(n)}(\sigma_0+\lambda)$ varies regularly at $0$ with index $\rho<0$\,. Then
\[
	\lim_{t\to\infty}\frac{1}{t}\log\nu(t,\infty)=\sigma_0\,.
\]
\end{Cor}

The following result is a generalization of the result of Nakagawa (cf. \cite[Theorem 3]{Na}).
\begin{Cor}\label{cor:tail}
	Let $f\in \cm$ with the representing measure $\nu$ and the abscissa of convergence $\sigma_0\in (-\infty,0]$ and assume that $\sigma_0$ is a pole of $f$. Then
	\[
		\lim_{t\to\infty}\frac{1}{t}\log\nu(t,\infty)=\sigma_0\,.
	\]	
\end{Cor}

Corollary \ref{cor:tail} has a simple  application in the probability theory. Let $X$ be a non-negative random variable defined on a probability space $(\Omega,\sig,\P)$ and define $\varphi(\lambda):= \E[e^{-\lambda X}]:= \int_\Omega e^{-\lambda X(\omega)} \P(d\omega)$. Then $\varphi$ is the Laplace transform of the law of the random variable $X$ and $\varphi\in \cm$ with the representing measure $\nu(dt)=\P(X\in dt)$\,. If the abscissa of convergence $\sigma_0$ of $\varphi$ satsfies $\sigma_0\in (-\infty,0]$ and $\varphi$ has a pole at $\sigma_0$, then it follows from Corollary \ref{cor:tail} that
\[
 \lim_{x\to\infty}\frac{1}{x}\log\P(X>x)=\sigma_0\,.
\]
This is the main result in \cite{Na}\,, where complex methods were used (as in the proof of Graham-Vaaler's theorem with minorant and majorant functions, see \cite{GV}).  Unlike approach in these articles, we use methods from real analysis. 

Results of this paper can treat more general examples. Namely, condition that $\phi$ has a pole at $\sigma_0$ can be relaxed to a more  general case of  regular variation at $\sigma_0$ with a negative index or, more generally, scaling properties given in (A-3). Furthermore, our results treat probability measures but also more general measures, e.g. L\' evy measures.

In the results presented so far it can be seen that only the exponential term dominates. Under additional assumptions it is possible to obtain finer asymptotical properties of the representing measure. First, let us introduce the following conditions for $f\in \cm$ with the abscissa of convergence $\sigma_0$:
\begin{itemize}
	      \item[{\bf (A-1)}] the representing measure  $\nu$ of $f$ has a
density with
respect to the Lebesgue measure, i.e. there exists a function 
\[
      \nu\colon (0,\infty)\rightarrow (0,\infty)\ \text{ such that }\ \nu(dt)=\nu(t)\,dt\,;
\]
		\item[{\bf (A-2)}] $\sigma_0\in (-\infty,0]$ and there exists $\eta\in \R$ such that $t\mapsto e^{-\sigma_0 t}t^\eta\nu(t)$ is monotone;
		\item[{\bf (A-3)}] there exist constants $c>0$, $0\leq \Lambda_1<\Lambda_2\leq \infty$, $n\in \N_0$ and $\gamma>0$ such
that
		\[
			\frac{(-1)^nf^{(n)}(\sigma_0+\lambda_2)}{(-1)^nf^{(n)}(\sigma_0+\lambda_1)}\leq c\,
\left(\frac{\lambda_2}{\lambda_1}\right)^{-\gamma}\
\text{ for all }\ \Lambda_1<\lambda_1\leq \lambda_2<\Lambda_2\,.
		\]	
	\end{itemize}
For example, (A-3) will hold in the case of regular variation with the help of Potter's theorem (see \cite[Theorem 1.5.6 (iii)]{BGT}).
To be more precise, it holds if $\lambda \mapsto (-1)^n f^{(n)}(\sigma_0+\lambda)$ varies regularly with index $\rho<0$ 
at the origin (take $\Lambda_1=0$ and $\Lambda_2<\infty$)
or
at infinity (take $\Lambda_1>0$ and $\Lambda_2=\infty$)\,.

\begin{Thm}\label{thm:main}
	Let $f\in \cm$ and assume that it satisfies (A-1). 

	(i) If (A-2)  holds, then
	there is a constant $c_1>0$ such that
	\[
		\nu(t)\leq c_1 (-1)^nf^{(n)}(\sigma_0+t^{-1})t^{-n-1}e^{\sigma_0 t}\ \text{ for all }\
t>
0\,.
	\]	

	(ii) If  (A-2) and (A-3) hold, then there exist constants $c_2>0$ and $\delta\in (0,1)$ such that 
	\[
		\nu(t)\geq c_2(-1)^nf^{(n)}(\sigma_0+t^{-1}) t^{-n-1}e^{\sigma_0 t}\ \text{ for all }\
t\in (\delta^{-1}\Lambda_2^{-1},\delta\Lambda_1^{-1})\,.
	\]
\end{Thm}

Recall that the notation $f(t)\asymp g(t), t\in I$ means that the quotient $f(t)/g(t)$ stays bounded from below and above for $t\in I$\,.
\begin{Cor}\label{cor:tail2}
Assume that  $f\in \cm$ satisfies (A-1)--(A-3) with $\Lambda_1=0$. Then there exists $t_0>0$ such that the density of the representing measure $\nu$ satisfies
\[
	\nu(t)\asymp (-1)^nf^{(n)}(\sigma_0+t^{-1}) t^{-n-1}e^{\sigma_0 t}\quad \text{ for all }\quad t\geq t_0\,.
\]
\end{Cor}
Compared to Theorem \ref{thm:tail}, the estimate in the last corollary  is more precise, but we had to assume more. 

Assumptions (A-1) and (A-2) seem to be rather technical. Nevertheless, they may be easily checked with the help of the following family of functions. We say that a $C^\infty$-function $\phi\colon (0,\infty)\rightarrow (0,\infty)$ is a Bernstein function if
$(-1)^{n+1}\phi^{(n)}(\lambda)\geq 0$ for all $\lambda>0$ and $n\in \N$\,.The class of Bernstein functions will be denoted by $\bbf$ and it is known that every $\phi\in \bbf$ can be uniquely represented in the following way (see \cite[Theorem 3.2]{SSV})
\begin{equation}\label{eq:bern_repr-intro}
 \phi(\lambda)=a+b\lambda+\int_{(0,\infty)}(1-e^{-\lambda t})\mu(dt)\,,
\end{equation}
where $a,b\geq 0$, $\mu$ is a measure on $(0,\infty)$ satisfying $\int_{(0,\infty)}(1\wedge t)\mu(dt)<\infty$\, usually called the L\' evy measure of $\phi$\,. There is a subclass of  \bbf\ that will play an important role known as complete Bernstein functions denoted by \cbf, which comprises of Bernstein functions $\phi\in \bbf$ whose L\' evy measure in the representation (\ref{eq:bern_repr-intro}) has a completely monotone density (with respect to the Lebesgue measure)\,. For example, $\lambda^\alpha \,\,(0\leq\alpha\leq 1)$ and $\log(1+\lambda)$ belong to $\cbf$\,. Taking derivative in (\ref{eq:bern_repr-intro}) we get 
\begin{equation}\label{eq:intro-tmp2}
	\phi'(\lambda)=\int_0^\infty e^{-\lambda t}t\mu(dt)\,.
\end{equation}
Hence, $\phi'\in \cm$ and it has the representing measure $\nu(dt)=t\mu(dt)$ by the uniqueness of the Laplace transform.   Also, composition of (complete) Bernstein functions stays (complete) Bernstein function. 

%
%

For $f\in \cm$ with the abscissa of convergence $\sigma_0$ we introduce the following conditions:
\begin{itemize}
	      \item[{\bf (B-1)}] the function $\lambda \mapsto \lambda f(\sigma_0+\lambda)$ is in $\bbf$;
	      \item[{\bf (B-2)}] $\sigma_0<0$ and there exists $a_0>0$ such that the function $\lambda\mapsto \int_0^\lambda f(\sigma_0+a+t)\,dt$ is in $\cbf$ for all $a\in (0,a_0)$\,.
\end{itemize}

\begin{Prop}\label{prop:Bs}
\begin{itemize}
	\item[(i)] If (B-1) holds, then (A-1) holds and $t\mapsto e^{-\sigma_0 t}\nu(t)$ is non-increasing.
	\item[(ii)] If (B-2) holds, then (A-1) holds and $t\mapsto t^{-1}e^{-\sigma_0 t}\nu(t)$ is completely monotone and, in particular, non-increasing. 
\end{itemize}
\end{Prop}

Let us illustrate our results by a few examples in which the measures are not explicitly known. 

\begin{Ex} (a)
	Let $\phi\in \bbf$ and assume that the abscissa of convergence of $f=\phi'\in \cm$ satisfies $\sigma_0\in (-\infty,0)$. Furthermore, assume that, for some $a_0>0$, 
	\[
		\lambda\mapsto \phi(\sigma_0+a+\lambda)-\phi(\sigma_0+a)\,\,\text{ is in }\,\, \cbf \,\,\text{ for all }\,\, a\in (0,a_0)\,.
	\]
	Then (B-2) holds, since
	\[
		\int_0^\lambda f(\sigma_0+a+t)\,dt=\int_{\sigma_0+a}^{\sigma_0+a+\lambda} \phi'(t)\,dt=\phi(\sigma_0+a+\lambda)-\phi(\sigma_0+a)\text{ is in }\,\,\cbf\,
	\]
	and this together with Proposition \ref{prop:Bs} implies (A-1) and (A-2)\,.
	Hence, if $f=\phi'$ satisfies (A-3), it follows from Theorem \ref{thm:main} and (\ref{eq:intro-tmp2}) that the L\' evy density of $\phi$ satisfies
	\[
		\mu(t)\asymp (-1)^n \phi^{(n+1)}(\sigma_0+t^{-1})t^{-n-2}e^{\sigma_0 t},\quad  t\in (\delta^{-1}\Lambda_2^{-1},\delta\Lambda_1^{-1})\,.
	\]
	(b) As a concrete example, let us consider $\phi(\lambda)=\log(1+\log(1+\lambda))$\,. In this case the L\' evy measure is not known, but we will obtain its behavior at $0$ and infinity. We have $\sigma_0=e^{-1}-1$ and for any $a>0$
	\[
		\phi(\sigma_0+a+\lambda)-\phi(\sigma_0+a)= 
		\log\left(1+\frac{\log\left(1+\frac{\lambda}{e^{-1}+a}\right)}{\log(1+ae)}\right)\in \cbf;
	\]
	hence (B-2) holds. Note that (A-3) holds with $n=0$, $\gamma<1$, $\Lambda_1=0$ and $\Lambda_2=\infty$, since
	\[
		f(\sigma_0+\lambda)=\frac{1}{(e^{-1}+\lambda)\log(1+e\lambda)}\asymp  \begin{cases}
			\frac{1}{\lambda} & 0<\lambda\leq 2\\
			\frac{1}{\lambda\log\lambda} & \lambda> 2\\
		\end{cases}
	\]
	and therefore
	\[
		\mu(t)\asymp \begin{cases}
			\frac{1}{t\log\frac{1}{t}} & 0<t<\tfrac{1}{2}\\
			\frac{e^{-(1-e^{-1})t}}{t} & t\geq \tfrac{1}{2}\,.
		\end{cases}
	\]
(c) 	Example (b) can be generalized by iterating logarithms. Define
\[
 \phi_1(\lambda):=\log(1+\lambda)\ \text{ and }\ \phi_{n+1}:=\phi_n\circ\phi_1\ \text{
for }\ n\in \N\,.
\]
Using the approach from (b) it follows that the L\' evy density $\mu_n(t)$ of $\phi_n$ satisfies
\[
 \mu_n(t)\asymp \begin{cases}
                 \frac{1}{t}\prod\limits_{k=1}^{n-1}\frac{1}{\phi_k(t^{-1})} & 0<t<\tfrac{1}{2}\\
		  \frac{e^{\phi_n^{-1}(-1)t}}{t} & t\geq \tfrac{1}{2}\,.
                \end{cases}
\]
Notice that $\lim\limits_{n\to\infty}\phi_n^{-1}(-1)=0$, meaning that, by iterating, the rate of exponential decay of the corresponding L\' evy measure becomes very close to zero. 
\end{Ex}

It is left to investigate the case when the abscissa of convergence is zero.
\begin{Ex}\label{ex:heavy1}
 Let $f(\lambda)=\lambda^{\gamma-1}$, where $\gamma\in (0,1)$.  Then the abscissa of convergence of $f$ is $\sigma_0=0$. In this case, exponential decay cannot be expected. In fact, the representing measure is explicitly known $\nu(dt)=\frac{dt}{\Gamma(1-\gamma)t^\gamma}$. 
 
 \end{Ex}
The following result explores such situations. Note that if the representing measure $\nu$ of $f\in \cm$ is finite, then $f(0)=\nu([0,\infty))$ is also finite\,. If $\lim\limits_{t\to\infty}\frac{\log\nu(t,\infty)}{\log t}\in (-\infty,0)$ we say that the tail of the measure $\nu$ is heavy. 
\begin{Thm}\label{thm:heavy}
Let $f\in \cm$ with the abscissa of convergence $\sigma_0=0$ and such that the representing measure $\nu$ is finite. Assume that there exist $\gamma>0$ and $n\in \N$ satisfying
\begin{equation}\label{eq:cond-htail1}
	\limsup_{\lambda\to0+}\frac{\log(f(0)-f(\lambda))}{\log\tfrac{1}{\lambda}}\leq -\gamma\,,\qquad\qquad \liminf_{\lambda\to0+}\frac{\log \left[(-1)^n f^{(n)}(\lambda)\right]}{\log\frac{1}{\lambda}}\geq -\gamma+n
\end{equation}
and
\begin{equation}\label{eq:cond-htail2}
	\limsup_{\lambda\to0+}\frac{f^{(n)}(2\lambda)}{f^{(n)}(\lambda)}<1\,.
\end{equation}
Then
\[
	\lim_{t\to\infty}\frac{\log\nu(t,\infty)}{\log t}=-\gamma\,.
\]
\end{Thm}

\begin{Ex}
 Let us consider $f\in \cm$ of the form $f(\lambda)=\alpha(\lambda)\lambda^\gamma\log\lambda+\beta(\lambda)$, where $\alpha$ and $\beta$ are analytic in the neighborhood of $0$, $\gamma\in (0,1]$, $\alpha(0)\not=0$ and such that the abscissa of convergence of $f$ is $\sigma_0=0$ and $f(0)=\beta(0)=1$. Then the conditions of Theorem \ref{thm:heavy} hold with $n=1$ for $\gamma<1$ and $n=2$ for $\gamma=1$, since 
 \[
  \lim_{\lambda\to0+}\frac{\log(f(0)-f(\lambda))}{\log\tfrac{1}{\lambda}}\leq \lim_{\lambda\to0+}\frac{\log\left(\frac{1-\beta(\lambda)}{\lambda^\gamma}-\alpha(\lambda)\log\lambda\right)}{\log\tfrac{1}{\lambda}}-\gamma=-\gamma
 \]
 and
 \begin{align*}
  \lim_{\lambda\to 0+}\frac{\log\left[-f'(\lambda)\right]}{\log\frac{1}{\lambda}}&=\lim_{\lambda\to 0+}\frac{\log\left(-\alpha(\lambda)\log\lambda-\alpha'(\lambda)\lambda\log\lambda+\alpha(\lambda)+\beta'(\lambda)\lambda^{1-\gamma}\right)}{\log\frac{1}{\lambda}}-\gamma+1\\&=-\gamma+1\,.
 \end{align*}
  To check (\ref{eq:cond-htail2}), for $\gamma<1$, 
  \begin{align*}
    \limsup_{\lambda\to 0+}\frac{f'(2\lambda)}{f'(\lambda)}&=\limsup_{\lambda\to 0+}\frac{\left(2\alpha'(2\lambda)\lambda+\alpha(2\lambda)2^{\gamma-1}\right)\frac{\log(2\lambda)}{\log\lambda}+\alpha(2\lambda)\frac{2^{\gamma-1}}{\log\lambda}}{\alpha'(\lambda)\lambda+\alpha(\lambda)+\frac{\alpha(\lambda)}{\log\lambda}}\\&=\frac{\alpha(0)2^{\gamma-1}}{\alpha(0)}=2^{\gamma-1}<1\,.
  \end{align*}
 In the case $\gamma=1$ we can similarly check that $\limsup_{\lambda\to 0+}\frac{f''(2\lambda)}{f''(\lambda)}=\frac{1}{2}<1\,.$
 Hence, the tail of the representing measure satisfies
 \[
  \lim_{t\to\infty}\frac{\log\nu(t,\infty)}{\log t}=-\gamma\,.
 \]
 Example \ref{ex:heavy1} is a special case of this example\,.

\end{Ex}


The structure of the paper is as follows. In Section \ref{sec:exp} main results concerning exponential decay are proved, while in Section \ref{sec:bf} these results are applied to the class of Bernstein functions. The result concerning heavy tails is proven in \ref{sec:heavy}. In Section \ref{sec:tauber} we prove Theorem \ref{thm:main}. Some applications of our results are given in Section \ref{sec:app}. 
In the first application we examine whether random sum of identically distributed heavy tailed random variables remains heavy tailed and determine the rate. The second application is in the theory of non-local equations. 
More precisely, we investigate exponential decay of fundamental solutions of some non-local equations such as
\[
	\log(1-\Delta)u+u=f\quad \text{ in }\quad \R^d\,,
\]
where $\Delta$ is the Laplacian in $\R^d$\,.
It turns out that the fundamental solution decays exponentially with rate $-\sqrt{1-e^{-1}}$ (see Example \ref{ex:log-fund})\,.

\section{Exponential decay}\label{sec:exp}

In this section the proof of Theorem  \ref{thm:tail} is given. The proof relies on the fact that the limit $\lim\limits_{t\to\infty}\frac{\log\nu(t,\infty)}{t}$ exists since the function $t\mapsto \frac{\log\nu(t,\infty)}{t}$ is monotone. 

\proof[Proof of Theorem \ref{thm:tail}]

By (\ref{eq:extension}), for any $\lambda\in (0,-\sigma_0)$ we have
\begin{align*}
(-1)^nf^{(n)}(\sigma_0+\lambda)&\geq \int_{\lambda^{-1}}^\infty e^{(-\sigma_0-\lambda )t}t^n \nu(dt)\\
&\geq e^{-\sigma_0\lambda^{-1}-1}\lambda^{-n}\nu(\lambda^{-1},\infty)
\end{align*}
implying
\begin{equation}\label{eq:upper-tail}
	\nu(t,\infty)\leq (-1)^nf^{(n)}(\sigma_0+t^{-1})e^{\sigma_0 t+1}t^{-n}\quad \text{ for }\quad t>(-\sigma_0)^{-1}\,.
\end{equation}
Hence,
\begin{align}
\lim_{t\to\infty} &\frac{1}{t}\log\nu(t,\infty)\leq \limsup_{t\to\infty}\left\{\frac{1}{t}\log\left[(-1)^nf^{(n)}(\sigma_0+t^{-1})\right]+\sigma_0+\frac{1-n\log{t}}{t}\right\}=\sigma_0\,.\label{eq:limsup}
\end{align}

To prove the equality in (\ref{eq:limsup}) in the last display, for any $\lambda>0$  we define the "tilted" measure $\nu_\lambda$ by
\[
	\nu_\lambda(dt)=e^{-\sigma_0t-\lambda t}t^n\nu(dt)\,.
\]	
It is a finite measure, since
\[
	\nu_\lambda(0,\infty)=\int_{(0,\infty)}e^{-\sigma_0 t-\lambda t}t^n\nu(dt)=(-1)^nf^{(n)}(\sigma_0+\lambda)<\infty\,.
\]
Furthermore, by Fubini theorem the following holds
\begin{align}
	\lambda \int_0^\infty e^{-\lambda t}\nu_\lambda(t,\infty)\,dt&=\int_{(0,\infty)}\int_0^s \lambda e^{-\lambda t}\,dt\,\nu_\lambda(ds)=\int_{(0,\infty)}(1-e^{-\lambda s})\,\nu_\lambda(ds)\nonumber\\
	&=\nu_\lambda(0,\infty)-\int_{(0,\infty)}e^{-\sigma_0 s-2\lambda s}s^n\nu(ds)\nonumber\\
	&=(-1)^nf^{(n)}(\sigma_0+\lambda)-(-1)^nf^{(n)}(\sigma_0+2\lambda)\,.\label{eq:tail-tmp3}
\end{align}
Let $\delta>0$ be chosen so that 
\begin{equation}\label{eq:tail-tmp4}
	e^{-\delta}>\limsup_{\lambda\to0+}\frac{(-1)^nf^{(n)}(\sigma_0+2\lambda)}{(-1)^nf^{(n)}(\sigma_0+\lambda)}\,.
\end{equation}
Now we split the integral 
\begin{align*}
	\lambda\int_0^\infty e^{-\lambda t}\nu_\lambda(t,\infty)\,dt&=\lambda \int_0^{\delta \lambda^{-1}}e^{-\lambda t}\nu_\lambda(t,\infty)\,dt+\lambda \int_{\delta \lambda^{-1}}^\infty e^{-\lambda t}\nu_\lambda(t,\infty)\,dt\\
	&\leq \nu_\lambda(0,\infty)\lambda \int_0^{\delta \lambda^{-1}}e^{-\lambda t}\,dt+\nu_\lambda(\delta\lambda^{-1},\infty)\lambda \int_{\delta \lambda^{-1}}^\infty e^{-\lambda t}\,dt\\
	&=(-1)^nf^{(n)}(\sigma_0+\lambda)(1-e^{-\delta})+\nu_\lambda(\delta\lambda^{-1},\infty)e^{-\delta}
\end{align*}
and use  (\ref{eq:tail-tmp3}) to conclude
\[
	e^{-\delta}(-1)^nf^{(n)}(\sigma_0+\lambda)-(-1)^nf^{(n)}(\sigma_0+2 \lambda)\leq e^{-\delta}\nu_\lambda(\delta\lambda^{-1},\infty)\,.
\]
Note that (\ref{eq:tail-tmp4}) and the last display imply
\begin{equation}\label{eq:tail-tmp5}
	\liminf_{\lambda\to0+}\frac{\nu_\lambda(\delta\lambda^{-1},\infty)}{(-1)^nf^{(n)}(\sigma_0+\lambda)}>0\,.
\end{equation}
Assume that $
	\lim\limits_{t\to\infty}\frac{1}{t}\log\nu(t,\infty)<\sigma_0\,.
$
Then there exist $t_0>1$ and $\varepsilon>0$ such that 
\begin{equation}\label{eq:tail-tmp10}
	\nu(t,\infty)\leq e^{(\sigma_0-\varepsilon)t}\quad \text{ for all }\quad t\geq t_0\,.
\end{equation}
Using integration by parts (or Fubini theorem) and (\ref{eq:tail-tmp10})  it follows that for some constant $c>0$
\begin{align*}
	\nu_\lambda(\delta\lambda^{-1},\infty)&
	=\nu(\delta\lambda^{-1},\infty)e^{-\sigma_0\delta\lambda^{-1}-\delta}(\delta\lambda^{-1})^n\\&\,\,\,\,\,\,+\int_{\delta\lambda^{-1}}^\infty \hspace{-.4cm}e^{-\sigma_0 t-\lambda t}t^{n-1}((-\sigma_0-\lambda)t+n)\nu(t,\infty)\,dt\\
							&\leq ce^{-\varepsilon \lambda^{-1}-\delta}+c(-\sigma_0+n)\int_{\delta\lambda^{-1}}^\infty e^{-\varepsilon t-\lambda t}t^{n-1}\,dt
							\quad \text{ for }\quad \lambda\in (0,\delta t_0^{-1})\,,
\end{align*}
yielding
$
	\liminf\limits_{\lambda\to 0+}\nu_\lambda(\delta\lambda^{-1},\infty )=0.
$
This contradicts (\ref{eq:tail-tmp5}), since $\lim\limits_{\lambda\to 0+}(-1)^nf^{(n)}(\sigma_0+\lambda)>0$\,. Hence, 
\[
 \lim\limits_{t\to\infty}\frac{1}{t}\log\nu(t,\infty)=\sigma_0\,.
\]
\qed


It is left to prove consequences of this theorem. 

\proof[Proof of Corollary \ref{cor:tail_reg-var}]
	Let $\ell\colon (0,\infty)\rightarrow  (0,\infty)$ be a slowly varying function so that
	\[
		(-1)^nf^{(n)}(\sigma_0+\lambda)=\lambda^\rho\ell(\lambda)\quad \text{ for all }\quad \lambda>0\,.
	\]
	Using the representation theorem for slowly varying functions (cf. $0$ - version of \cite[Theorem 1.3.1]{BGT}), there exist $a>0$ and measurable functions $c,\varepsilon :(0,a)\rightarrow \R$ so that $\lim\limits_{\lambda\to0+} c(\lambda)>0$, $\lim\limits_{\lambda\to0+} \varepsilon(\lambda)=0$ and
	\[
		\ell(\lambda)=c(\lambda)\exp{\left\{\int_\lambda^a\frac{\varepsilon(u)}{u}\,du\right\}}\quad \text{ for all }\quad \lambda \in (0,a)\,.
	\]
	Hence, if $\varepsilon_0>0$ is chosen so that $|\varepsilon(\lambda)|\leq \varepsilon_0$ for all $\lambda\in (0,a)$ we have 
	\begin{align*}
		0&\leq \liminf_{\lambda\to 0+} \lambda\log\left[(-1)^nf^{(n)}(\sigma_0+\lambda)\right]\\&\leq \limsup_{\lambda\to 0+}\left[\rho\lambda\log\lambda +\lambda\log c(\lambda) +\varepsilon_0\lambda \log\frac{a}{\lambda}\right]=0\,.
	\end{align*}
	The other condition follows directly from the definition of the regular variation:
	\[
		\limsup_{\lambda\to0+}\frac{(-1)^nf^{(n)}(\sigma_0+2\lambda)}{(-1)^nf^{(n)}(\sigma_0+\lambda)}=2^{\rho}<1\,.
	\]
\qed

\proof[Proof of Corollary \ref{cor:tail}]
	Let 
	\[
		f(\lambda)=\sum_{k=1}^m\frac{a_k}{(\lambda-\sigma_0)^k}+f_0(\lambda)
	\]
	be the Laurent series expansion of $f$ around $\sigma_0$, where $f_0$ is analytic in $\sigma_0$\, and $m\in \N$ is the order of the pole. Note that $a_m>0$, since $\sigma_0$ is the pole of order $m$ and $\lim\limits_{\lambda\to 0+} f(\sigma_0+\lambda)=\infty$. Then we can easily check the function $h\colon (0,\infty)\rightarrow (0,\infty)$ defined by $h(\lambda)=f(\sigma_0+\lambda)$ varies regularly at $0$ with index $-m$, since
	\begin{align*}
		\lim_{\lambda\to0+}\frac{h(\lambda x)}{h(\lambda)}&=\lim_{\lambda\to 0+}\frac{a_mx^{-m}+\lambda \sum\limits_{k=1}^{m-1}a_kx^{-k}\lambda^{m-1-k}+\lambda^mf_0(\lambda x)}{a_m+\lambda \sum\limits_{k=1}^{m-1}a_k\lambda^{m-1-k}+\lambda^mf_0(\lambda )}\\&=x^{-m}\quad \text{ for any }\quad x>0\,.
	\end{align*}
	Therefore,  it is possible to apply Corollary \ref{cor:tail_reg-var} with $n=0$\,.
\qed

\section{Tauberian theorem}\label{sec:tauber}
We start with a result that relates condition (A-3) to the second condition in (\ref{eq:cond-tail}) in some cases. 
\begin{Lem}\label{lem:equiv}
Let $f\colon (0,\infty)\rightarrow (0,\infty)$ be a non-increasing function. The following claims are equivalent: 
\begin{itemize}
	\item[(i)] there exists $b>1$ such that $\limsup\limits_{\lambda\to0+}\frac{f(b\lambda)}{f(\lambda)}<1\,; $
	\item[(ii)] there exist $c>0$, $\gamma>0$ and $\lambda_0>0$ such that 
	\[
		\frac{f(\lambda_2)}{f(\lambda_1)}\leq c\left(\frac{\lambda_2}{\lambda_1}\right)^{-\gamma}\text{ for all } 0<\lambda_1\leq \lambda_2\leq \lambda_0 \text{ and } x\geq 1\,.
	\]
\end{itemize}
\end{Lem}
\proof
Assume that (ii) holds. Choosing $b>c^{1/\gamma}\vee 1$ we see that 
\[
    \limsup_{\lambda\to 0+}\frac{f(b\lambda)}{f(\lambda)}\leq c b^{-\gamma}<1.
\]
Assume now that (i) holds. There exists $\lambda_0>0$ such that $\kappa:=\sup\limits_{\lambda\in (0,\lambda_0)}\frac{f(b\lambda)}{f(\lambda)}\in (0,1)$. Let $0<\lambda_1\leq \lambda_2\leq \lambda_0$. There exists $n\in \N$ such that $b^{n-1}\leq \frac{\lambda_2}{\lambda_1}<b^n$\,. Since $f$ is non-increasing, we obtain
	\[
		\frac{f(\lambda_2 )}{f(\lambda_1)}\leq \frac{f(b^{n-1} \lambda_1)}{f(\lambda_1)}\leq \kappa^{n-1}=\kappa^{-1}\left(b^{n}\right)^{\frac{\log\kappa}{\log b}}\leq \kappa^{-1} \left(\frac{\lambda_2}{\lambda_1}\right)^{\frac{\log\kappa}{\log b}}\,.
	\]
	Here we have used that $\gamma:=-\frac{\log\kappa}{\log b}>0$\,. 
\qed

\proof[Proof of Theorem \ref{thm:main}]
(i) First we assume that the function in (A-2) is non-increasing. If $\eta>0$, then the function $t\mapsto e^{-\sigma_0 t}\nu(t)$ is also non-increasing; thus we may assume that $\eta\leq 0$. 
  Using  (\ref{eq:extension}) and (A-2), for $\lambda>0$, we get 
\begin{align*}
 (-1)^nf^{(n)}(\sigma_0+\lambda)& \geq \int\limits_0^{\lambda^{-1}}e^{-\sigma_0 t}e^{-\lambda
t}t^n\nu(t)\,dt \geq\lambda^{-\eta}
\nu(\lambda^{-1})e^{-\sigma_0 \lambda^{-1}}\int\limits_0^{\lambda^{-1}}e^{-\lambda t}t^{n-\eta}\,dt &
\\
&= \lambda^{-n-1}e^{-\sigma_0 \lambda^{-1}}\nu(\lambda^{-1})\int_0^1t^{n-\eta}e^{-t}\,dt\,.
\end{align*}
This gives the upper bound; take  $t>0$ and set $\lambda=t^{-1}$ to deduce from the previous
display that
\begin{equation}\label{eq:sec2_9}
 \nu(t)\leq c_1t^{-n-1}e^{\sigma_0 t}(-1)^nf^{(n)}(\sigma_0+t^{-1}),
\end{equation}
where $c_1=\left(\int_0^1t^{n-\eta}e^{-t}\,dt\right)^{-1}$\,. If the function in (A-2) is non-decreasing, we may proceed similarly by using the estimate 
\[
 (-1)^nf^{(n)}(\sigma_0+\lambda) \geq \int\limits_{\lambda^{-1}}^\infty e^{-\sigma_0 t}e^{-\lambda t}t^n\nu(t)\,dt
\]
to obtain the same  bound as in (\ref{eq:sec2_9}) with $c_1=\left(\int_1^\infty t^{n-\eta}e^{-t}\,dt\right)^{-1}$\,.
\\
(ii) Assume that the function in (A-2) is non-increasing. Let $\delta\in (0,1)$ and $\lambda\in (\Lambda_1,\Lambda_2)$. Then
{\allowdisplaybreaks
\begin{align*}
 \int\limits_{\delta\lambda^{-1}}^\infty e^{-\sigma_0 t}e^{-\lambda
t}t^n\nu(t)\,dt&=(-1)^nf^{(n)}(\sigma_0+\lambda)-\int\limits_0^{\delta\lambda^{-1}}e^{-\sigma_0
t}e^{-\lambda t}t^n\nu(t)\,dt \\
&\geq (-1)^nf^{(n)}(\sigma_0+\lambda)-c_1\int\limits_0^{\delta\lambda^{-1}}e^{-\lambda t}
t^{-1}(-1)^nf^{(n)}(\sigma_0+t^{-1})\,dt\\
&\geq (-1)^nf^{(n)}(\sigma_0+\lambda)-c_1
(-1)^nf^{(n)}(\sigma_0+\lambda)\lambda^\gamma
\int\limits_0^{\delta\lambda^{-1}}e^{-\lambda t}t^{\gamma-1}\,dt \\
&\geq (-1)^nf^{(n)}(\sigma_0+\lambda)-c_1\gamma^{-1}
\delta^{\gamma}(-1)^nf^{(n)}(\sigma_0+\lambda)\,.
\end{align*}}
Choosing $\delta\in (0,1)$ small enough so that $1-c_1\gamma^{-1}
\delta^{\gamma}\geq \frac{1}{2}$ one obtains
{\allowdisplaybreaks
\begin{align*}
 \tfrac{1}{2}(-1)^nf^{(n)}(\sigma_0+\lambda) & \leq  \int\limits_{\delta\lambda^{-1}}^\infty
e^{-\sigma_0 t}e^{-\lambda
t}t^n\nu(t)\,dt\leq (\delta\lambda^{-1})^{\eta}\nu(\delta \lambda^{-1})e^{-\sigma_0
\lambda^{-1}}\int\limits_{\delta\lambda^{-1}}^\infty e^{-\lambda t}t^{n-\eta}\,dt \\
&=(\delta\lambda^{-1})^{-n-1}\nu(\delta\lambda^{-1})e^{-\sigma_0\delta\lambda^{-1}}\int\limits_1^\infty e^{-\delta^{-1}t}t^{n-\eta}\,dt\,.
\end{align*}}
Let $t\in (\delta^{-1}\Lambda_2^{-1},\delta\Lambda_1^{-1})\subset (\Lambda_2^{-1},\Lambda_1^{-1})$. Then $\lambda=\delta t^{-1}\in (\Lambda_1,\Lambda_2)$ and thus 
the last display and (A-3) imply
\begin{align*}
 \nu(t)&\geq c_2t^{-n-1}(-1)^n f^{(n)}(\sigma_0+\delta t^{-1})e^{\sigma_0 t} 
\geq c_3 t^{-n-1}(-1)^nf^{(n)}(\sigma_0)+t^{-1})e^{\sigma_0 t}.
\end{align*}
The case when the function in (A-2) is non-decreasing can be proven similarly. 
\qed

\section{Heavy tails}\label{sec:heavy}
In this section we prove theorem concerning the behavior of the tail when the abscissa of convergence is zero. 

\proof[Proof of Theorem \ref{thm:heavy}]
Since $f(0)=\nu(0,\infty)<\infty$ we can perform a similar calculation as in (\ref{eq:tail-tmp3}) to obtain
\begin{align}
	\lambda\int_0^\infty e^{-\lambda t}\nu(t,\infty)\,dt&=\int_{(0,\infty)}\int_0^s\lambda e^{-\lambda t}\,dt\nu(ds)
								     =\int_{(0,\infty)}(1-e^{-\lambda s})\nu(ds)\nonumber\\
								    &=\nu(0,\infty)-f(\lambda)=f(0)-f(\lambda)\,.\label{eq:heavy-tmp0}
\end{align}
Hence, 
\begin{equation*}\label{eq:heavy-tmp50}
f(0)-f(\lambda)\geq \lambda\int_0^{\lambda^{-1}}e^{-\lambda t}\nu(t,\infty)\,dt\geq e^{-1}\nu(\lambda^{-1},\infty)\,,
\end{equation*}
which gives
\begin{equation}\label{eq:heavy-upper}
	\nu(t,\infty)\leq e(f(0)-f(t^{-1}))\quad \text{ for }\quad t>0\,. 
\end{equation}
This implies 
\[
	\lim_{t\to\infty}\frac{\log\nu(t,\infty)}{\log t}\leq \limsup_{t\to\infty }\frac{\log(f(0)-f(t^{-1}))}{\log t}=-\gamma\,.
\]
For any $\lambda>0$ we define a measure $\nu_\lambda(dt)=e^{-\lambda t}t^n\,\nu(dt)$. By Fubini theorem we have
\begin{align}
\lambda \int_0^\infty e^{-\lambda t}\nu_\lambda(t,\infty)\,dt&=\int_0^\infty e^{-\lambda s}s^n\int_0^s \lambda e^{-\lambda t}\,dt\nu(ds)\nonumber\\
&=\int_0^\infty(e^{-\lambda s}-e^{-2\lambda s})s^n\nu(ds)=(-1)^nf^{(n)}(\lambda)-(-1)^nf^{(n)}(2\lambda)\,.\label{eq:heavy-tmp60}
\end{align}
On the other hand, if we choose $\delta>0$ so that 
\begin{equation}\label{eq:heavy-tmp62}
e^{-\delta}>\limsup_{\lambda\to 0+}\frac{(-1)^nf^{(n)}(2\lambda)}{(-1)^nf^{(n)}(\lambda)},
\end{equation}
we obtain
\begin{align}
 \lambda \int_0^\infty e^{-\lambda t}\nu_\lambda(t,\infty)\,dt&= \lambda \int_0^{\delta\lambda^{-1}}e^{-\lambda t}\nu_\lambda(t,\infty)\,dt+ \lambda \int_{\delta\lambda^{-1}}^\infty e^{-\lambda t}\nu_\lambda(t,\infty)\,dt\nonumber\\
 &\leq \nu_\lambda(0,\infty)\lambda \int_0^{\delta\lambda^{-1}} e^{-\lambda t}\,dt+\nu_\lambda(\delta\lambda^{-1},\infty)\int_{\delta\lambda^{-1}}^\infty \lambda e^{-\lambda t}\,dt\nonumber\\
 &=(-1)^nf^{(n)}(\lambda)(1-e^{-\delta})+\nu_\lambda(\delta\lambda^{-1},\infty)e^{-\delta}\,.\label{eq:heavy-tmp61}
\end{align}
Then it follows from (\ref{eq:heavy-tmp60}), (\ref{eq:heavy-tmp61}) and the choice of $\delta$ in (\ref{eq:heavy-tmp62}) that 
\begin{equation}\label{eq:heavu-tmp66}
  \liminf_{\lambda\to0+}\frac{\nu_\lambda(\delta\lambda^{-1},\infty)}{(-1)^nf^{(n)}(\lambda)}\geq e^{\delta}\left(e^{-\delta}-\limsup_{\lambda\to 0+}\frac{(-1)^nf^{(n)}(2\lambda)}{(-1)^nf^{(n)}(\lambda)}\right)>0\,.
\end{equation}

Assume that $\lim\limits_{t\to\infty}\frac{\log\nu(t,\infty)}{\log t}<-\gamma$. Then there exist $\varepsilon>0$ and $t_0>0$ such that 
$ \nu(t,\infty)\leq t^{-\gamma-\varepsilon}$ for $t\geq t_0$. Using integration by parts, for $\lambda>0$ large enough, 
\begin{align*}
 \nu_\lambda(\delta_\lambda^{-1},\infty)&=\int\limits_{\delta\lambda^{-1}}^\infty e^{-\lambda t}t^n\nu(dt)
 = e^{-\delta}(\delta\lambda^{-1})^n\nu(\delta\lambda^{-1},\infty)+\int\limits_{\delta\lambda^{-1}}^\infty (n-\lambda t) t^{n-1}\nu(t,\infty)\,dt\\
 &\leq  e^{-\delta} (\delta\lambda^{-1})^{n-\gamma-\varepsilon}+n\int\limits_{\delta\lambda^{-1}}^{n\lambda^{-1}}e^{-\lambda t}t^{n-1-\gamma-\varepsilon}\,dt\leq c_1 \lambda^{-n+\gamma+\varepsilon}\,.
\end{align*}
Hence, by (\ref{eq:cond-htail1}) it follows that 
\[
  \liminf_{\lambda\to0+}\frac{\nu_\lambda(\delta\lambda^{-1},\infty)}{(-1)^nf^{(n)}(\lambda)}\leq \liminf_{\lambda\to0+}\frac{c_1\lambda^{-n+\gamma+\varepsilon}}{\lambda^{\gamma-n+\varepsilon/2}}=0
\]
contradicting (\ref{eq:heavu-tmp66}). Therefore, $\lim\limits_{t\to\infty}\frac{\log\nu(t,\infty)}{\log t}=-\gamma$\,.
\qed


\section{Bernstein functions}\label{sec:bf}

In this section we apply proven results to Bernstein functions. Let us recall that a Bernstein function $\phi\in \bbf$ has the following representation
\begin{equation}\label{eq:bern_repr}
 \phi(\lambda)=a+b\lambda+\int_{(0,\infty)}(1-e^{-\lambda t})\mu(dt)\,,
\end{equation}
where $a,b\geq 0$ and  $\mu$ is the L\' evy measure. 

\begin{Prop}\label{prop:bf-tail}
 Let $\phi\in \bbf$ and assume that the abscissa of convergence $\sigma_0$ of $f=\phi'$ satisfies $\sigma_0\in (-\infty,0]$\,. If $f$ satisfies the assumptions of Theorem \ref{thm:tail}, Corollary \ref{cor:tail_reg-var} or Corollary \ref{cor:tail}, then
 \[
  \lim_{t\to\infty}\frac{\log \mu(t,\infty)}{t}=\sigma_0\,,
 \]
 where $\mu$ is the L\' evy measure of $\phi$ in the representation (\ref{eq:bern_repr-intro})\,.
\end{Prop}
\proof
Using (\ref{eq:intro-tmp2}) and conditions of the proposition it follows that the representing measure of $\phi'$ is $\nu(dt)=t\mu(dt)$; hence we obtain by
Theorem \ref{thm:tail}, Corollary \ref{cor:tail_reg-var} or Corollary \ref{cor:tail} the following 
\begin{equation}\label{eq:tmp-tail-bf1}
 \lim_{t\to\infty}\frac{1}{t}\log\left(\int_t^\infty s\mu(ds)\right)=\sigma_0\,.
\end{equation}
From this we deduce directly  that 
\[
 \sigma_0\geq \limsup_{t\to\infty}\frac{1}{t}\log(t\mu(t,\infty))=\lim_{t\to\infty}\frac{1}{t}\log\mu(t,\infty)\,.
\]
Assume that the strict inequality holds. Then there exist $t_0>0$ and $\varepsilon\in (0,-\sigma_0)$ so that 
\[
 \mu(t,\infty)\leq e^{(\sigma_0-\varepsilon)t}\quad \text{ for }\quad t\geq t_0\,.
\]
This would imply
\begin{align*}
 \int_t^\infty s\mu(ds)&= t\mu(t,\infty)+ \int_t^\infty \mu(s,\infty)\,ds\\
 &\leq te^{(\sigma_0-\varepsilon)t}+(-\sigma_0-\varepsilon)^{-1}e^{(\sigma_0+\varepsilon)t}\quad \text{ for }\quad t\geq t_0\,, 
\end{align*}
and consequently $ \lim\limits_{t\to\infty}\frac{1}{t}\log\left(\int_t^\infty s\mu(ds)\right)\leq \sigma_0-\varepsilon$ contradicting (\ref{eq:tmp-tail-bf1})\,.
\qed

\proof[Proof of Proposition \ref{prop:Bs}] (i) Assume that (B-1) holds. 
Then there exists a measure $\mu$ on $(0,\infty)$ and $a,b\geq 0$ such that 
\[
	\lambda f(\lambda+\sigma_0)=a+b\lambda+\int_{(0,\infty)}(1-e^{-\lambda t})\mu(dt)
\]
Since $\frac{1-e^{-\lambda t}}{\lambda}\leq 1\wedge t$ for $\lambda\geq 1$ we can use dominated convergence theorem to obtain 
\[
    \frac{1}{\lambda}\lim_{\lambda\to\infty}\int_{(0,\infty)}(1-e^{-\lambda t})\mu(dt)=0.
\]
Therefore, since $\sigma_0>-\infty$, 
\[
	b=\lim_{\lambda\to \infty}\frac{\lambda f(\sigma_0+\lambda)}{\lambda}=\lim_{\lambda \to\infty }\int_{(0,\infty)}e^{-(\sigma_0+\lambda)t}\nu(dt)=0\,.
\]
It follows that 
\[
	f(\sigma_0+\lambda)=\frac{a}{\lambda}+\int_0^\infty e^{-\lambda t}\mu(t,\infty)\,dt=\int_0^\infty e^{-\lambda t}(a+\mu(t,\infty))\,dt,
\]
which, by the uniqueness of the Laplace transform, implies that the representing measure $\nu$ of $f$ has a density (hence satisfies (A-1)) given by
\[
	\nu(t)=e^{\sigma_0 t}(a+\mu(t,\infty))
\]
and it satisfies (A-2) with $\eta=0$\,. 
\\
(ii) If (B-2) holds, then for any $a\in (0,a_0)$ there exists $\mu_a\in \cm$ so that
\begin{align*}
 \int\limits_0^\infty (1-e^{-\lambda t})\mu_a(t)\,dt&=\int\limits_0^\lambda f(\sigma_0+a+t)\,dt=\int\limits_0^\lambda \int\limits_0^\infty e^{-\sigma_0 s-as-st}\nu(ds)\,dt\\
 &=\int\limits_0^\infty (1-e^{-\lambda s})s^{-1}e^{-as}e^{-\sigma_0 s}\nu(ds)\,.
\end{align*}
Taking derivative and using the uniqueness of the Laplace transform it follows that $\nu(dt)=\nu(t)\,dt$ with
\[
 \nu(t)=te^{at} e^{\sigma_0 t}\mu_a(t),\quad t>0\,.
\]
Using the fact that $\cm$ is closed under pointwise limits (see \cite[Corollary 1.6]{SSV}), it follows that 
\[
 t^{-1}e^{-\sigma_0 t}\nu(t)=\lim_{a\to 0+}t^{-1}e^{-\sigma_0 t}e^{-at}\nu(t)=\lim_{a\to 0+}\mu_a(t)\in \cm\,.
\]
\qed

Let us make a short excursus and link Bernstein functions to a class of stochastic processes called subordinators. A stochastic process $S=\{S_t:t\geq 0\}$ defined on a probability space $(\Omega,\sig,\P)$ is called a subordinator if it takes values in $[0,\infty)$, it has stationary and independent increments (i.e., for any $0\leq t_1\leq t_2\leq \ldots\leq t_n,\ n\in \N$ the random variables $S_{t_1},S_{t_2}-S_{t_1},\ldots,S_{t_n}-S_{t_{n-1}}$ are independent and identically distributed) and $\P(\{\omega\in \Omega\colon t\mapsto S_t(\omega) \text{ is right-continuous with left limits}\})=1$\,. In other words, $S$ is a L\' evy process taking values in $[0,\infty)$. It turns out that in this case we can calculate the Laplace transform of $S_t$ and it is of the form
\[
	\E[e^{-\lambda S_t}]=e^{-t\phi(\lambda)},\quad \lambda>0\,,
\]
where $\phi$ has a representation (\ref{eq:bern_repr}), i.e. $\phi\in \bbf$ (cf. \cite[Chapter 3]{Be})\, and in this context it is usually called the Laplace exponent of $S$. The L\' evy measure $\mu$ represents intensity of jumps of the process $S$. 

\begin{Ex}\label{ex:rel}
 Let $S=\{S_t\colon t\geq 0\}$ be a relativistic $\alpha/2$-stable subordinator, i.e. $\phi(\lambda)=(\lambda+m^{2/\alpha})^{\alpha/2}-m$ for some $\alpha\in (0,2)$ and $m>0$\,. Then the Laplace transform of $S_t$ is $f(\lambda)=e^{-t\phi(\lambda)}$ with the abscissa of convergence  $\sigma_0=-m^{2/\alpha}$\,. Since $h(\lambda)=-tf'(\sigma_0+\lambda)=te^{-\lambda^{\alpha/2}-m}\lambda^{\alpha/2-1}$ varies regularly at $0$ with index $\alpha/2-1<0$, it follows from Corollary \ref{cor:tail_reg-var} that 
 \[
  \lim_{r\to\infty}\frac{\log\P(S_t> r)}{r}=-m^{2/\alpha}\,.
 \]
\end{Ex}

\section{Applications}\label{sec:app}
\subsection{Random sums}
In this subsection we investigate whether random sum of heavy tailed random variables remains heavy tailed. Let $Y=\{Y_i:i\geq 1\}$ be a sequence of independent and identicaly distributed random variables taking values in $[0,\infty)$ and assume that function   
\[
f\in \cm \quad \text{ defined by }\quad f(\lambda)=\E[e^{-\lambda Y_1}] \ \text{ for }\ \lambda>0
\]
satisfies conditions (\ref{eq:cond-htail1}) and (\ref{eq:cond-htail2}) for some $\gamma>0$ and $n=1$\,.
Let $N$ be random variable taking values in the set $\N$ that is independent of the sequence $Y$ and denote by $p_m:=\P(N=m)$, $m\in \N$ its law. 
The random sum is defined by $X=\sum\limits_{m=1}^N Y_m\,.$
It is easy to see that the Laplace transform of the random variable $X$ is given by
\[
g(\lambda)=\sum_{m=1}^\infty f(\lambda)^m p_m\,.
\]
Let us investigate how the tail of $X$ behaves in two situations. 

(a)\quad  If we assume that $\E N=\sum\limits_{m=1}^\infty mp_m<\infty$, then it is easy to see that $g$ satisfies conditions of Theorem \ref{thm:heavy}\,. Indeed, since $\frac{1-f(\lambda)^m}{1-f(\lambda)}\leq m$, we get 
\begin{align*}
\limsup_{\lambda\to0+}\frac{\log(g(0)-g(\lambda))}{\log\frac{1}{\lambda}}&=\limsup_{\lambda \to 0+}\frac{\log\left(\sum\limits_{m=1}^\infty \frac{1-f(\lambda)^m}{1-f(\lambda)}p_m\right)+\log(f(0)-f(\lambda))}{\log\frac{1}{\lambda}}\leq -\gamma\,.
\end{align*}
Also, since $N$ has finite mean, by the dominated convergence theorem,
\begin{align*}
 \liminf_{\lambda\to 0+}\frac{\log\left[-g'(\lambda)\right]}{\log\frac{1}{\lambda}}&=\liminf_{\lambda\to 0+}\frac{\log\left[-f'(\lambda)\right]+\log\left(\sum_{m=1}^\infty f(\lambda)^{m-1}mp_m\right)}{\log\frac{1}{\lambda}}\\
 &=\liminf_{\lambda\to 0+}\frac{\log\left[-f'(\lambda)\right]+\log\E N}{\log\frac{1}{\lambda}}\geq -\gamma+1
\end{align*}
and
\begin{align*}
 \limsup_{\lambda\to 0+}\frac{g'(2\lambda)}{g'(\lambda)}&=\limsup_{\lambda\to 0+}\frac{f'(2\lambda)}{f'(\lambda)}\cdot \frac{\sum\limits_{m=1}^\infty f(\lambda)^{m-1}m p_m}{\sum\limits_{m=1}^\infty f(2\lambda)^{m-1}m p_m}=\limsup_{\lambda\to 0+}\frac{f'(2\lambda)}{f'(\lambda)}\cdot \frac{\E N}{\E N}<1\,.
\end{align*}

Hence, Theorem \ref{thm:heavy} implies that $X$ is heavy tailed with the same rate $-\gamma$ as $Y_1$\,. 

(b)\quad Assume that $p_m=cm^{-2}$, $m\in \N$,  where $c>0$ is the normalizing constant. Note that in this case $N$ does not have finite mean. In this special case it follows that
\[
	g(\lambda)=-c\int_0^{f(\lambda)}\frac{\log(1-t)}{t}\,dt\quad \text{ for }\quad \lambda>0\,.
\]
Since
\[
	g(0)-g(\lambda)\leq \frac{c}{f(\lambda)}\int_{f(\lambda)}^1\left(-\log(1-t)\right)\,dt=(1-f(\lambda))\left(1-\log(1-f(\lambda))\right),
\]
it follows that for any $\varepsilon>0$ there is a constant $c>0$ such that 
\begin{align*}
	\limsup_{\lambda\to 0+}\frac{\log(g(0)-g(\lambda))}{\log\frac{1}{\lambda}}&=\limsup_{\lambda\to 0+}\frac{\log(1-f(\lambda))+\log\left(1-\log(1-f(\lambda))\right)}{\log\frac{1}{\lambda}}\\
	&\leq \limsup_{\lambda\to 0+}\frac{\log(1-f(\lambda))+\log\left(1+c(1-f(\lambda)^{-\varepsilon}\right)}{\log\frac{1}{\lambda}}\\
	&\leq -\gamma(1-\varepsilon)
\end{align*}
implying that $\limsup\limits_{\lambda\to 0+}\frac{\log(g(0)-g(\lambda))}{\log\frac{1}{\lambda}}\leq -\gamma$\,.

Furthermore, since $t\mapsto -c\frac{\log(1-t)}{t}$ is non-decreasing on $(0,1)$, it follows that $-c\frac{\log(1-f(\lambda))}{f(\lambda)}\geq c$ and therefore
\begin{align*}
 \liminf_{\lambda\to 0+}\frac{\log(-g'(\lambda))}{\log\frac{1}{\lambda}}&=\liminf_{\lambda\to 0+}\frac{\log\left[-c\frac{\log(1-f(\lambda))}{f(\lambda)}(-f'(\lambda))\right]}{\log\frac{1}{\lambda}}\\&\geq \liminf_{\lambda\to 0+}\frac{\log\left[-cf'(\lambda)\right]}{\log\frac{1}{\lambda}}\geq -\gamma+1\,.
\end{align*}

To check the last condition of Theorem \ref{thm:heavy}, first we note that 
\[
 \frac{g'(2\lambda)}{g'(\lambda)}=\frac{\log(1-f(2\lambda))}{\log(1-f(\lambda))}\cdot\frac{f(\lambda)}{f(2\lambda)}\cdot\frac{f'(2\lambda)}{f'(\lambda)}\,.
\]
Now, since $1-f(2\lambda)=\int\limits_0^{2\lambda}(-f'(t)\,dt)\geq -2\lambda f'(2\lambda)$,
\[
 \frac{\log(1-f(2\lambda))}{\log(1-f(\lambda))}\leq \frac{\log(2\lambda)+\log(-f'(\lambda))}{\log(1-f(\lambda))},
\]
hence, by (\ref{eq:cond-htail1}) and (\ref{eq:cond-htail2}), it follows that 
\begin{align*}
 \limsup_{\lambda\to 0+}\frac{g'(2\lambda)}{g'(\lambda)}&\leq \limsup_{\lambda\to 0+}\frac{\frac{\log(2\lambda)}{\log\frac{1}{\lambda}}+\frac{\log(-f'(\lambda))}{\log\frac{1}{\lambda}}}{\frac{\log(1-f(\lambda))}{\log\frac{1}{\lambda}}}\limsup_{\lambda\to 0+}\frac{f'(2\lambda)}{f'(\lambda)}\\&\leq\frac{-1-\gamma+1}{-\gamma}\limsup_{\lambda\to 0+}\frac{f'(2\lambda)}{f'(\lambda)}<1\,.
\end{align*}
Therefore,  $g$ satisfies the assumptions of Theorem \ref{thm:heavy}; hence $X$ is heavy tailed and \[\lim\limits_{x\to+\infty}\frac{\log\P(X>x)}{\log{x}}=-\gamma\,.\]

\subsection{Non-local equations}
As a 
second 
application of our main results we consider rate of decay of solutions of the 
equation 
\begin{equation}\label{eq:sec5_1}
 \phi(-\Delta)u=f\ \text{ in }\ \R^d\,,
\end{equation}
where $\Delta=\frac{\partial^2}{\partial x_1^2}+\frac{\partial^2}{\partial x_2^2}+\ldots+\frac{\partial^2}{\partial x_d^2}$ is the Laplacian in $\R^d$ and $\phi\in \bbf$.

The operator $\phi(-\Delta)$ should be understood in terms of the Fourier transform. More precisely, the
Fourier transform of $f\in L^1(\R^d)$ is defined by
\[ \wh{f}(\xi):=\int\limits_{\R^d}e^{i\xi\cdot x}f(x)\,dx,\ \xi\in \R^d
\]
and we use the same notation for the extension of the Fourier transform from $L^1(\R^d)\cap L^2(\R^d)$
to a unitary operator on $L^2(\R^d)$ (see \cite[Theorem 8.29]{Fo}). 

Now the operator $\phi(-\Delta)$ in equation  (\ref{eq:sec5_1}) is understood as a pseudo-differential operator:
\begin{equation}\label{eq:nonloc}
 [\phi(-\Delta)u]^{\ \wh{}}(\xi):=\phi(|\xi|^2) \wh{u}(\xi),\ \xi\in \R^d\,,
\end{equation}
for 
$
u\in D(\phi(-\Delta)):=\{u\in
L^2(\R^d)\colon \phi(|\xi|^2) \wh{u}(\xi)\ \text{ is in }\ L^2(\R^d)\}\,.
$
In this subsection we will investigate decay of the fundamental solution of the equation (\ref{eq:sec5_1}). This is a function $K\colon
\R^d\setminus \{0\}\rightarrow \R$ defined by
\begin{equation}\label{eq:fund_soln}
 K(x)=\int\limits_{(0,\infty)} p(t,x)\nu(dt),\ x\in
\R^d\setminus \{0\}\,,
\end{equation}
where $\nu$ is the representing measure of the function $\frac{1}{\phi}\in \cm$ and  $p\colon (0,\infty)\times \R^d \rightarrow \R$ is the Gauss-Weierstrass kernel
defined by
\[
p(t,x)=(4\pi t)^{-d/2} e^{-\frac{|x|^2}{4t}}\,,\quad t>0,\,x\in \R^d\,.
\]

\begin{Rem}\label{rem:gw} Since
 \[
\int\limits_{\R^d}p(t,x)\,dx=1\ \ \text{ and }\ \ \wh{p(t,\cdot)}(\xi)=e^{-t|\xi|^2}\,.
 \]
 it is easy to check that $u\colon \R^d\rightarrow \R$ defined by 
 \begin{equation}\label{eq:solution}
 u(x):=(K\star f)(x):=\int\limits_{\R^d}K(x-y)f(y)\,dy\ \text{ for }\ x\in \R^d
\end{equation}
solves the equation (\ref{eq:sec5_1}):
 \[
 [\phi(-\Delta)u]^{\ \wh{}}(\xi)=\phi(|\xi|^2)\wh{K}(\xi)\wh{f}(\xi)=\wh{f}(\xi)\,.
\]
\end{Rem}
Probabilistacally, $K$ can be understood as a potential (or Green function of the whole space $\R^d$) of the subordinate Brownian motion, where the Laplace exponent of the corresponding subordinator is $\phi$. 

The following proposition gives a link between exponential decay of fundamental solution and the tail of the representing measure. 
\begin{Prop}\label{prop:fund}
 Let $\phi\in \bbf$ and let $f=\frac{1}{\phi}\in \cm$ with the abscissa of convergence $\sigma_0\in (-\infty,0)$ and the representing measue $\nu$. If $\lim\limits_{t\to\infty}\frac{\log\nu(t,\infty)}{t}=\sigma_0$, then
 \[
  \lim_{|x|\to\infty}\frac{\log K(x)}{|x|}=-\sqrt{-\sigma_0}\,.
 \]
\end{Prop}
\proof
Using integration by parts and change of variable we get from (\ref{eq:fund_soln}) the following
\begin{equation}\label{eq:non-loc-eq1}
 K(x)\leq\int\limits_0^\infty (4\pi t)^{-d/2}e^{-\frac{|x|^2}{4t}}\nu(t,\infty)\frac{|x|^2}{4t^2}\,dt
 =\pi^{-d/2}|x|^{-d}\int\limits_0^\infty e^{-t} t^{d/2} \nu\left(\frac{|x|^2}{4t},\infty\right)\,dt
\end{equation}


and hence, for any $\varepsilon\in (0,-\sigma_0)$, the assumptions imply
\begin{align}
 \limsup_{|x|\to\infty}\frac{\log K(x)}{|x|}&\leq  \limsup_{|x|\to\infty}\frac{\log\left[\pi^{-d/2}|x|^{-d}\int\limits_0^\infty e^{-t+(\sigma_0+\varepsilon)\frac{|x|^2}{4t}}\,t^{d/2}\,dt\right]}{|x|}\nonumber\\
 &=\limsup_{|x|\to\infty}\frac{\log\left[\int\limits_0^\infty e^{-(\sqrt{t}-\sqrt{-\sigma_0-\varepsilon}\frac{|x|}{\sqrt{t}})^2}\,t^{d/2}\,dt\, e^{-\sqrt{-\sigma_0-\varepsilon}|x|}\right]}{|x|}=-\sqrt{-\sigma_0-\varepsilon}\,.\label{eq:non-loc-upper}
\end{align}
The last equality can be seen if we perform change of variable $y=-\sqrt{t}+\sqrt{-\sigma_0-\varepsilon}\frac{|x|}{\sqrt{t}}$ in the integral. Then we get  $t=\left(\frac{-y+\sqrt{y^2+4\sqrt{-\sigma_0-\varepsilon}|x|}}{2}\right)^2$ and, hence, 
\begin{align*}
	\int\limits_0^\infty e^{-\left(\sqrt{t}-\sqrt{-\sigma_0-\varepsilon}\frac{|x|}{\sqrt{t}}\right)^2}&t^{d/2}\,dt=\int\limits_{-\infty}^\infty e^{-y^2}\left(\frac{-y+\sqrt{y^2+4\sqrt{-\sigma_0-\varepsilon}|x|}}{2}\right)^{d+1}\hspace{-0.5cm}\frac{dy}{\sqrt{y^2+4\sqrt{-\sigma_0-\varepsilon}|x|}},
\end{align*}
implying, by the dominated convergence theorem, the following
\[
 \lim_{|x|\to\infty}|x|^{-d/2}\int\limits_0^\infty e^{-\left(\sqrt{t}-\sqrt{-\sigma_0-\varepsilon}\frac{|x|}{\sqrt{t}}\right)^2}t^{d/2}\,dt=\int\limits_{-\infty}^\infty e^{-y^2}\frac{(-\sigma_0-\varepsilon)^{\frac{d+1}{4}}}{2\sqrt{\sqrt{-\sigma_0-\varepsilon}}}\,dy
 =\frac{\sqrt{\pi}}{2}(-\sigma_0-\varepsilon)^{\frac{d}{4}}\,.
\]
Since $\varepsilon\in (0,-\sigma_0)$ was arbitrary we obtain from (\ref{eq:non-loc-upper})
\[
\limsup_{|x|\to\infty}\frac{\log K(x)}{|x|}\leq -\sqrt{-\sigma_0}\,.
\]	
To obtain the lower bound, integration by parts yields 
\begin{align*}
K(x)&=\int\limits_0^\infty (4\pi t)^{-\frac{d}{2}-1}e^{-\frac{|x|^2}{4t}}\left(\frac{|x|^2}{4t}-\frac{d}{2}\right)\nu(t,\infty)\,dt\\
&\geq \frac{d}{2}\int\limits_0^{\frac{|x|^2}{4d}} (4\pi t)^{-\frac{d}{2}-1}e^{-\frac{|x|^2}{4t}}\nu(t,\infty)\,dt-\frac{d}{2}\int\limits_{\frac{|x|^2}{4d}}^\infty (4\pi t)^{-\frac{d}{2}-1}e^{-\frac{|x|^2}{4t}}\nu(t,\infty)\,dt
\end{align*}
Let $\varepsilon\in (0,-\sigma_0)$. Then there exist constants $c_1,c_2>0$ such that  $c_1e^{(\sigma_0-\varepsilon)t}\leq \nu(t,\infty)\leq c_2e^{(\sigma_0+\varepsilon)t}$ for $t>0$. 
Hence,
\begin{equation}\label{eq:K-tmp2}
	K(x)\geq c_3 e^{-\sqrt{-\sigma_0+\varepsilon}|x|}\int\limits_0^{\frac{|x|^2}{2d}}t^{-\frac{d}{2}-1}e^{-\left(\sqrt{-\sigma_0+\varepsilon}\sqrt{t}-\frac{|x|}{2\sqrt{t}}\right)^2}\,dt-c_4 |x|^{-d-2}e^{\frac{\sigma_0+\varepsilon}{2d}|x|^2}\,.
\end{equation}
By the change of variable $s=\frac{|x|^2}{4t}$ in the integral in the last display we obtain 
\[
	I(x):=\int\limits_0^{\frac{|x|^2}{2d}}t^{-\frac{d}{2}-1}e^{-\left(\sqrt{-\sigma_0+\varepsilon}\sqrt{t}-\frac{|x|}{2\sqrt{t}}\right)^2}\,dt=2^d|x|^{-d}\int\limits_{\sqrt{d}}^\infty t^{\frac{d}{2}-1}e^{-\left(\sqrt{-\sigma_0+\varepsilon}\frac{|x|}{\sqrt{t}}-\sqrt{t}\right)^2}\,dt
\]
and then using the change of variable $y=-\sqrt{t}+\sqrt{-\sigma_0+\varepsilon}\frac{|x|}{\sqrt{t}}$ similar to the one in the first part of the proof it follows that for $|x|$ large enough
\[
	I(x)\geq \int\limits_{-\infty}^{0}e^{-y^2}\left(\frac{-y+\sqrt{y^2+4\sqrt{-\sigma_0+\varepsilon}|x|}}{2}\right)^{d}\frac{dy}{\sqrt{y^2+4\sqrt{-\sigma_0+\varepsilon}|x|}}\,.
\]
Hence, by the dominated convergence theorem, 
\[
	\liminf_{|x|\to\infty} |x|^{-\frac{d-1}{2}}I(x)\geq \int\limits_{-\infty}^0 e^{-y^2}\frac{\left(-\sigma_0+\varepsilon\right)^\frac{d}{4}}{\sqrt{\sqrt{-\sigma_0+\varepsilon}}}\,dy=\frac{\sqrt{\pi}}{2}\left(-\sigma_0+\varepsilon\right)^\frac{d-1}{4}\,.
\]
Using the last display in (\ref{eq:K-tmp2}), we conclude that, for $|x|$ large enough,
\[
	K(x)\geq c_6 |x|^{\frac{d-1}{2}} e^{-\sqrt{-\sigma_0+\varepsilon}|x|},
\]
hence
\[
	\liminf_{|x|\to\infty}\frac{\log K(x)}{|x|}\geq -\sqrt{-\sigma_0+\varepsilon}\,.
\]
Now the lower bound follows since $\varepsilon\in (0,\sigma_0)$ was arbitrary.

\qed
\begin{Ex}\label{ex:log-fund} Consider the following equation
	\[
		\log(1-\Delta)u+\beta u=f\,,
	\]
	where $\beta>0$\,.
	Here $\phi(\lambda)=\log(1+\lambda)+\beta\in \bbf$ and the abscissa of convergence of $f=1/\phi$ is $\sigma_0=e^{-\beta}-1$. Since $f(\sigma_0+\lambda)=\frac{1}{\log(1+e^\beta\lambda)}$ varies regularly at 0 with index $-1$, it follows  from Corollary \ref{cor:tail_reg-var} that the tail of the representing measure of $f$ satisfies $\lim_{t\to\infty}\frac{\log\nu(t,\infty)}{t}=e^{-\beta}-1$ and thus by Proposition \ref{prop:fund} the fundamental solution $K$ decays exponentially with rate $-\sqrt{1-e^{-\beta}}$\,.	
\end{Ex}
\begin{Ex}
	For equations connected with the potentials of  the generalized relativistic Schr\"{o}dinger operators (see Example \ref{ex:rel})
	\[
		(m^{2/\alpha}-\Delta)^{\alpha/2}u-m u+\beta u=f,
	\]
	where $\alpha\in (0,2)$ and $0<\beta<m$,  we can easily show as in the previous example that the fundamental solution decays exponentially with rate $-\sqrt{m^{2/\alpha}-(m-\beta)^{2/\alpha}}$\,. Indeed, since the abscissa of convergence of $f(\lambda)=\frac{1}{(\lambda+m^{2/\alpha})^{\alpha/2}-m+\beta}\in \cm$ is $\sigma_0=(m-\beta)^{2/\alpha}-m^{2/\alpha}$ and the function $f(\sigma_0+\lambda)=\frac{1}{((m-\beta)^{2/\alpha}+\lambda)^{\alpha/2}-(m-\beta)}$ varies regularly at $0$ with index $-1$, Corollary \ref{cor:tail_reg-var} and Proposition \ref{prop:fund} apply. 
\end{Ex}

\begin{center}
Acknowledgement\\

The author dedicates this article to Professors Z. Krajina and J. Paladino.
\end{center}


\providecommand{\bysame}{\leavevmode\hbox to3em{\hrulefill}\thinspace}
\providecommand{\MR}{\relax\ifhmode\unskip\space\fi MR }
\providecommand{\MRhref}[2]{%
  \href{http://www.ams.org/mathscinet-getitem?mr=#1}{#2}
}
\providecommand{\href}[2]{#2}

\end{document}